%
%
%
\newif\ifsect\newif\iffinal
\secttrue\finaltrue
\def\smallsect #1. #2\par{\bigbreak\noindent{\bf #1.}\enspace{\bf #2}\par
    \global\parano=#1\global\eqnumbo=1\global\thmno=1
    \nobreak\smallskip\nobreak\noindent\message{#2}}
\def\thm #1: #2{\medbreak\noindent{\bf #1:}\if(#2\thmp\else\thmn#2\fi}
\def\thmp #1) { (#1)\thmn{}}
\def\thmn#1#2\par{\enspace{\sl #1#2}\par
        \ifdim\lastskip<\medskipamount \removelastskip\penalty 55\medskip\fi}
\def\square{{\msam\char"03}}
\def\qedn{\thinspace\null\nobreak\hfill\square\par\medbreak}
\def\pf{\ifdim\lastskip<\smallskipamount \removelastskip\smallskip\fi
        \noindent{\sl Proof\/}:\enspace}
\def\itm#1{\par\indent\llap{\rm #1\enspace}\ignorespaces}

\def\bar#1{\overline{#1}}
\newcount\parano
\newcount\eqnumbo
\newcount\thmno
\newcount\versiono
\def\neweqt#1$${\xdef #1{(\number\parano.\number\eqnumbo)}
    \eqno #1$$
    \iffinal\else\rsimb#1\fi
    \global \advance \eqnumbo by 1}
\def\newthmt#1 #2: #3{\xdef #2{\number\parano.\number\thmno}
    \global \advance \thmno by 1
    \medbreak\noindent
    \iffinal\else\lsimb#2\fi
    {\bf #1 #2:}\if(#3\thmp\else\thmn#3\fi}
\def\neweqf#1$${\xdef #1{(\number\eqnumbo)}
    \eqno #1$$
    \iffinal\else\rlap{$\smash{\hbox{\hfilneg\string#1\hfilneg}}$}\fi
    \global \advance \eqnumbo by 1}
\def\newthmf#1 #2: #3{\xdef #2{\number\thmno}
    \global \advance \thmno by 1
    \medbreak\noindent
    \iffinal\else\llap{$\smash{\hbox{\hfilneg\string#1\hfilneg}}$}\fi
    {\bf #1 #2:}\if(#3\thmp\else\thmn#3\fi}
\def\inizia{\ifsect\let\neweq=\neweqt\else\let\neweq=\neweqf\fi
\ifsect\let\newthm=\newthmt\else\let\newthm=\newthmf\fi}
\def\bititolo{\empty}
\gdef\begin #1 #2\par{\xdef\titolo{#2}
\ifsect\let\neweq=\neweqt\else\let\neweq=\neweqf\fi
\ifsect\let\newthm=\newthmt\else\let\newthm=\newthmf\fi
\centerline{\titlefont\titolo}
\if\bititolo\empty\else\medskip\centerline{\titlefont\bititolo}
\xdef\titolo{\titolo\ \bititolo}\fi
\bigskip
\centerline{\bigfont \autore}
\if\istituto!\else\bigskip
\centerline{\istituto}
\centerline{\indirizzo}
\centerline{\email}\fi
\medskip
\centerline{#1~\anno}
\bigskip\bigskip
\ifsect\else\global\thmno=1\global\eqnumbo=1\fi}
\def\istituto{!}
\def\anno{2004}
\font\titlefont=cmssbx10 scaled \magstep1
\font\bigfont=cmr12
\font\eightrm=cmr8
\font\sc=cmcsc10
\font\bbr=msbm10
\font\sbbr=msbm7
\font\ssbbr=msbm5
\font\msam=msam10

\nopagenumbers
\binoppenalty=10000
\relpenalty=10000
\newfam\amsfam
\textfont\amsfam=\bbr \scriptfont\amsfam=\sbbr \scriptscriptfont\amsfam=\ssbbr
\let\de=\partial

\def\Hol{\mathop{\rm Hol}\nolimits}

\def\id{\mathop{\rm id}\nolimits}
\mathchardef\void="083F
\def\C{{\mathchoice{\hbox{\bbr C}}{\hbox{\bbr C}}{\hbox{\sbbr C}}
{\hbox{\sbbr C}}}}

\newcount\notitle
\notitle=1
\headline={\ifodd\pageno\rhead\else\lhead\fi}
\def\rhead{\ifnum\pageno=\notitle\iffinal\hfill\else\hfill\tt Version
\the\versiono; \the\day/\the\month/\the\year\fi\else\hfill\eightrm\titolo\hfill
\folio\fi}
\def\lhead{\ifnum\pageno=\notitle\hfill\else\eightrm\folio\hfill\autore\hfill
\fi}
\newbox\bibliobox
\def\setref #1{\setbox\bibliobox=\hbox{[#1]\enspace}
    \parindent=\wd\bibliobox}
\def\biblap#1{\noindent\hang\rlap{[#1]\enspace}\indent\ignorespaces}
\def\art#1 #2: #3! #4! #5 #6 #7-#8 \par{\biblap{#1}#2: {\sl #3\/}.
    #4 {\bf #5} (#6)\if.#7\else, \hbox{#7--#8}\fi.\par\smallskip}
\def\book#1 #2: #3! #4 \par{\biblap{#1}#2: {\bf #3.} #4.\par\smallskip}
\def\coll#1 #2: #3! #4! #5 \par{\biblap{#1}#2: {\sl #3\/}. In {\bf #4,}
#5.\par\smallskip}
\def\pre#1 #2: #3! #4! #5 \par{\biblap{#1}#2: {\sl #3\/}. #4, #5.\par\smallskip}
\versiono=6
\def\autore{Marco Abate${}^*$ and Filippo Bracci\footnote{${}^*$}{\rm Both authors
have been partially supported by Progetto MIUR di Rilevante Interesse Nazionale {\it
Propriet\`a geometriche delle variet\`a reali e complesse.}}}
%
%
%
%
\begin {October} Ritt's theorem and the Heins map in hyperbolic complex manifolds

{\narrower

{\sc Abstract.} Let $X$ be a Kobayashi hyperbolic complex manifold, and assume that $X$ does
not contain compact complex submanifolds of positive dimension (e.g., $X$ Stein). We shall
prove the following generalization of Ritt's theorem: every holomorphic self-map $f\colon
X\to X$ such that $f(X)$ is relatively compact in $X$ has a unique fixed point $\tau(f)\in
X$, which is attracting. Furthermore, we shall prove that $\tau(f)$ depends holomorphically
on~$f$ in a suitable sense, generalizing results by Heins, Joseph-Kwack and the second
author.

}

\smallsect 0. Introduction

\def\autore{Marco Abate and Filippo Bracci}The classical Wolff-Denjoy theorem (see, e.g.,
[A2, Theorem~1.3.9]) says that the sequence of iterates of a holomorphic self-map $f$ of the
unit disk $\Delta\subset\C$, except when $f$ is an elliptic automorphism of~$\Delta$ or the
identity, converges uniformly on compact subsets to a point $\tau(f)\in\bar{\Delta}$, the
{\sl Wolff point} of~$f$. Furthermore, if $\tau(f)\in\Delta$ then it is the unique fixed
point of~$f$; and if $\tau(f)\in\de\Delta$ then it is still morally fixed, in the sense that
$f(\zeta)$ tends to~$\tau(f)$ when $\zeta$ tends to~$\tau(f)$ non-tangentially.

In 1941, Heins [H] proved that the map
$\tau\colon\Hol(\Delta,\Delta)\setminus\{\id\}\to\bar\Delta$, associating to every elliptic
automorphism its fixed point and to any other map its Wolff point, is continuous. More than
half a century later, using the first author's version~(see [A1]) of the Wolff-Denjoy theorem
for strongly convex domains in~$\C^n$, Joseph and Kwack [JK] extended Heins' result to
strongly convex domains. 

In 2002, the second author started investigating further regularity properties of the Heins
map. If $D$ is a bounded domain in~$\C^n$, then $\Hol(D,D)$ is a subset of the complex
Banach space~$H^\infty(D)^n$ of $n$-uples of bounded holomorphic functions defined on~$D$;
so one may ask whether the Heins map, when defined, is holomorphic on some suitable open
subset of~$\Hol(D,D)$. And indeed, in [B] the second author proved that, when $D$ is
strongly convex, the Heins map is well-defined and holomorphic on~$\Hol_c(D,D)$, the open
subset of holomorphic self-maps of~$D$ whose image is relatively compact in~$D$.

The aim of this paper is to prove a similar result for the space $\Hol_c(X,X)$ of the
holomorphic self-maps of a Kobayashi hyperbolic Stein manifold whose image is relatively
compact in~$X$. First of all, we shall generalize the classical Ritt's
theorem, proving (Theorem~1.1) that every $f\in\Hol_c(X,X)$ admits a unique fixed point
$\tau(f)\in X$; therefore the Heins map $f\mapsto\tau(f)$ is well-defined and continuous
(Lemma~2.1).

To study further regularity properties of the Heins map, one apparently needs a complex
structure on~$\Hol_c(X,X)$. Unfortunately, we do not know whether such a structure exists in
general; so we shall instead prove (Theorem~2.3) that the Heins map is holomorphic when
restricted to any holomorphic family inside~$\Hol_c(X,X)$, a fact equivalent to $\tau$ being
holomorphic with respect to any sensible complex structure on~$\Hol_c(X,X)$. For instance,
we obtain (Corollary~2.4) that the Heins map is holomorphic on~$\Hol_c(D,D)$ for any bounded
domain~$D$ in~$\C^n$.

Finally, the first author would like to thank prof. Weiping Yin and the Capital Normal 
University of Beijing for the warm hospitality he enjoyed during his stay in China.

\smallsect 1. Ritt's theorem

Let $X$ be a complex manifold. We shall denote by $\Hol_c(X,X)$ the space of holomorphic
self-maps $f\colon X\to X$ of~$X$ such that $f(X)$ is relatively compact in $X$.

In 1920, Ritt~[R] proved that if $X$ is a non-compact Riemann surface then every
$f\in\Hol_c(X,X)$ has a unique fixed point $z_0\in X$. Furthermore, this fixed point is {\sl
attractive} in the sense that the sequence~$\{f^k\}$ of iterates of~$f$ converges, uniformly
on compact subsets, to the constant map~$z_0$. This theorem has been generalized to bounded
domains in~$\C^n$ by Wavre~[W]; see also Herv\'e~[He, p.~83]. Arguing as in~[A2,
Corollary~2.1.32] we shall now prove a far-reaching generalization of Ritt's theorem:

\newthm Theorem \Ritt: Let $X$ be a hyperbolic manifold with no compact complex submanifolds
of positive dimension. Then every $f\in\Hol_c(X,X)$ has a unique fixed point~$z_0\in X$.
Furthermore, the sequence of iterates of~$f$ converges, uniformly on compact subsets, to the
constant map~$z_0$.

\pf Since $X$ is hyperbolic, by [A3] the space $\Hol(X,X)$ of holomorphic self-maps of~$X$ is
relatively compact in the space $C^0(X,X^*)$ of continuous maps of~$X$ into the one-point
compactification~$X^*=X\cup\{\infty\}$, endowed with the compact-open topology. If
$f\in\Hol_c(X,X)$, this implies that the sequence of iterates of~$f$ is relatively compact
in~$\Hol(X,X)$, because~$f(X)\subset\subset X$.

Let then $\{f^{k_\nu}\}$ be a subsequence of $\{f^k\}$ converging to~$h_0\in\Hol(X,X)$. We can
also assume that $p_\nu=k_{\nu+1}-k_\nu$ and $q_\nu=p_\nu-k_\nu$ tend to~$+\infty$ as
$\nu\to+\infty$, and that there are $\rho_0$,~$g_0\in\Hol(X,X)$ such that
$f^{p_\nu}\to\rho_0$ and $f^{q_\nu}\to g_0$ in~$\Hol(X,X)$. Then it is easy to see that
$$
h_0\circ\rho_0=h_0=\rho_0\circ h_0\qquad\hbox{and}\qquad g_0\circ h_0=\rho_0=h_0\circ g_0,
$$
and so
$$
\rho_0^2=\rho_0\circ\rho_0=g_0\circ h_0\circ\rho_0=g_0\circ h_0=\rho_0.
$$
Thus $\rho_0$ is a holomorphic retraction, whose image is contained in the closure of~$f(X)$,
which is compact. This means (see Rossi~[Ro] and Cartan~[C]) that $\rho_0(X)$ is a compact
connected complex submanifold of~$X$, i.e., a point~$z_0\in X$. Therefore $\rho_0\equiv z_0$
and $z_0$ is a fixed point of~$f$, since $f$ clearly commutes with~$\rho_0$.

We are left to proving that $f^k\to z_0$, which implies in
particular that~$z_0$ is the only fixed point of~$f$. Since
$\{f^k\}$ is relatively compact in~$\Hol(X,X)$, it suffices to
show that~$z_0$ is the unique limit point of any converging
subsequence of~$\{f^k\}$. So let $\{f^{k_\mu}\}$ be a subsequence
converging toward a map $h\in\Hol(X,X)$. Arguing as before we find
a holomorphic retraction~$\rho\in\Hol(X,X)$ such that $h=\rho\circ
h$. Furthermore, $\rho$ must again be constant; but since it is
obtained as a limit of a subsequence of iterates of~$f$, it must
commute with~$\rho_0$, and this is possible if and only if
$\rho\equiv z_0$. But then $h=\rho\circ h\equiv z_0$ too, and we
are done.\qedn

In particular this theorem holds for hyperbolic Stein manifolds, because a Stein manifold has
no compact complex submanifolds of positive dimension.
\medbreak
{\it Remark 1.1:} If $f^k\to z_0$, then the spectral radius of~$df_{z_0}$ is strictly less
than one. Indeed, if $df_{z_0}$ had an eigenvalue~$\lambda\in\C$ with $|\lambda|\ge 1$, then
$d(f^k)_{z_0}$ would have $\lambda^k$ as eigenvalue, and $\lambda^k\not\to 0$ whereas
$d(f^k)_{z_0}\to O$.

\smallsect 2. The Heins map

Let $X$ be a hyperbolic manifold with no compact complex submanifolds of positive dimension.
The {\sl Heins map} of $X$ is the map $\tau\colon\Hol_c(X,X)\to X$ that associates to any
$f\in\Hol_c(X,X)$ its unique fixed point~$\tau(f)\in X$, whose existence is proved in
Theorem~\Ritt.

The first observation is that the Heins map is continuous:

\newthm Lemma \Heins: Let $X$ be a hyperbolic manifold with no compact complex submanifolds
of positive dimension. Then the Heins map $\tau\colon\Hol_c(X,X)\to X$ is continuous.

\pf Let $\{f_k\}\subset\Hol_c(X,X)$ be a sequence converging toward a map~$f\in\Hol_c(X,X)$;
we must show that $\tau(f_k)\to\tau(f)\in X$.

First of all, we claim that the set $\{\tau(f_k)\}$ is relatively compact in~$X$. Assume that
this is not true; then, up to passing to a subsequence, we can assume that the sequence
$\{\tau(f_k)\}$ eventually leaves any compact subset of~$X$. Now, the set $f(X)$ is
relatively compact in~$X$; we can then find an open set~$D$ in~$X$ such that
$$
f(X)\subset\subset D\subset\subset X.
$$
We have $\tau(f_k)\notin\bar{D}$ eventually; therefore for $k$ large enough we can
find~$R_k>0$ such that
$$
\bar{B\bigl(\tau(f_k),R_k\bigr)}\cap D=\void\qquad\hbox{and}\qquad
\bar{B\bigl(\tau(f_k),R_k\bigr)}\cap\de D\ne\void,
$$
where $B(z,R)$ is the ball of center~$z\in X$ and radius~$R>0$
with respect to the Kobayashi distance of~$X$. Choose
$z_k\in\bar{B\bigl(\tau(f_k),R_k\bigr)}\cap\de D$ for every $k$
large enough; since $\de D$ is compact, up to a subsequence we can
assume that $z_k\to z_0\in\de D$. In particular, then,
$f_k(z_k)\to f(z_0)\in f(X)\subset D$. But, on the other hand, we
have $f_k(z_k)\in\bar{ B\bigl(\tau(f_k),R_k\bigr)}\subset
X\setminus D$ for all $k$ large enough, because $\tau(f_k)$ is
fixed by~$f_k$ and the Kobayashi distance is contracted by
holomorphic maps; therefore $f(z_0)\in X\setminus D$,
contradiction.

So $\{\tau(f_k)\}$ is relatively compact in~$X$; to prove that $\tau(f_k)\to\tau(f)$ it
suffices to show that $\tau(f)$ is the unique limit point of the sequence~$\{\tau(f_k)\}$.
But indeed if $\tau(f_{k_\nu})\to x\in X$ we have
$$
f(x)=\lim_{\nu\to+\infty} f_{k_\nu}\bigl(\tau(f_{k_\nu})\bigr)=
\lim_{\nu\to+\infty} \tau(f_{k_\nu})=x;
$$
but $\tau(f)$ is the only fixed point of~$f$, and we are done.\qedn

As stated in the introduction, our aim is to prove that the Heins map is holomorphic in a
suitable sense. Since we do not know how to define a holomorphic structure on~$\Hol_c(X,X)$
for general manifolds, we shall prove another result which is equivalent to the holomorphy
of~$\tau$ in any reasonable setting (see for instance Corollary~2.4 below). We shall need the
following lemma:

\newthm Lemma \Abate: Let $P\subset\C^n$ be a polydisk centered in~$p_0\in\C^n$, and $h\colon
P\to\C^n$ a holomorphic map. Then there is a holomorphic map $A\colon P\to M(n,\C)$, where
$M(n,\C)$ is the space of $n\times n$ complex matrices, satisfying the following properties:
\smallskip
\itm{(i)} $h(z)-h(p_0)=A(z)\cdot(z-p_0)$ for all $z\in P$;
\itm{(ii)} $A(p_0)=dh_{p_0}$;
\itm{(iii)} for every polydisk $P_1\subset\subset P$ centered at~$p_0$ there is a constant
$C(P_1)>0$ such that $\|A\|_{P_1}\le C(P_1)\|h\|_{P_1}$.

\pf We can write
$$
h(z)-h(p_0)=\int_0^1{\de\over\de t} h\bigl(z_0+t(z-p_0)\bigr)\, dt=
\sum_{j=1}^n (z^j-p_0^j)\int_0^1 {\de h\over\de z^j}\bigl(z_0+t(z-p_0)\bigr)\, dt.
$$
Therefore taking
$$
A^i_j(z)=\int_0^1 {\de h^i\over\de z^j}\bigl(z_0+t(z-p_0)\bigr)\, dt
$$
the matrix $A=(A^i_j)$ clearly satisfies (i) and (ii), and (iii) follows from the Cauchy
estimates.\qedn

\newthm Theorem \Bracci: Let $X$ be a hyperbolic manifold with no compact complex submanifolds
of positive dimension, $Y$ another complex manifold, and $F\colon Y\times X\to X$ a
holomorphic map so that $f_y=F(y,\cdot)\in\Hol_c(X,X)$ for every $y\in Y$. Then the
map~$\tau_F\colon Y\to X$ given by $\tau_F(y)=\tau(f_y)$ is holomorphic. Furthermore, for
every~$y_0\in Y$ the differential of~$\tau_F$ at~$y_0$ is given by
$$
d(\tau_F)_{y_0}=\bigl(\id-d(f_{y_0})_{\tau(f_{y_0})}\bigr)^{-1}\circ dF_{(y_0,\tau(f_{y_0}))}
(\cdot,O).
$$

Notice that, by Remark~1.1, $\id-d(f_{y_0})_{\tau(f_{y_0})}$ is invertible.

\pf Without loss of generality, we can assume that $Y$ is a ball
$B^m\subset\C^m$ centered at~$y_0$. Set $p_0=\tau(f_{y_0})$, and
let $P_0\subset X$ be the domain of a polydisk chart centered
at~$p_0$. Since $f_{y_0}(p_0)=p_0$, we can find a polydisk
$P_1\subset\subset P_0$ centered at~$p_0$ such that
$f_{y_0}(P_1)\subset\subset P_0$. Furthermore, by Lemma~\Heins\
there is also a $\delta>0$ such that $\|y-y_0\|<\delta$ implies
$\tau(f_y)\in P_1$ and $f_y(P_1)\subset\subset P_0$. This means
that as soon as $y$ is close enough to~$y_0$ we can work
inside~$P_0$ and assume, without loss of generality, that $X$ is
contained in some~$\C^n$.

Write $p_y=\tau(f_y)\in P_1$, and define $h_y\colon\bar{P_1}\to\C^n$ by $h_y=f_y-f_{y_0}$.
We have
$$
p_y-p_0=f_{y_0}(p_y)-f_{y_0}(p_0)+h_y(p_y);
$$
therefore Lemma~\Abate\ applied to $f_{y_0}$ yields a
matrix~$A(y)$, depending continuously on~$y$ by Lemma \Heins, such
that
$$
p_y-p_0=A(y)\cdot(p_y-p_0)+h_y(p_y).
$$
Since $A(y)\to d(f_{y_0})_{p_0}$ as $y\to y_0$, for $y$ close to~$y_0$ $\id-A(y)$ is
invertible, and so we can write
$$
p_y-p_0=\bigl(\id-A(y)\bigr)^{-1}\cdot h_y(p_y).
\neweq\equno
$$

Now, we have
$$
dF_{(y_0,\tau(f_{y_0}))}(\cdot, O)=\hbox{Jac}_y\bigl(f_y(p_0)\bigr)(y_0),
$$
where $\hbox{Jac}_y$ is the Jacobian matrix computed with respect to the~$y$ variables; in
particular,
$$
h_y(p_0)-dF_{(y_0,\tau(f_{y_0}))}(y-y_0, O)=o(\|y-y_0\|).
$$
This means that to show that $\tau_F$ is holomorphic and $d\tau_F$ has the claimed expression
it suffices to show that
$$
\lim_{y\to y_0}{\bigl\|\tau_F(y)-\tau_F(y_0)-\bigl(\id-d(f_{y_0})_{p_0}\bigr)^{-1}\cdot
h_y(p_0)\bigr\|\over\|y-y_0\|}=0,
$$
which is equivalent to proving that
$$
\lim_{y\to y_0}{\bigl\|\bigl(\id-d(f_{y_0})_{p_0}\bigr)\cdot(p_y-p_0)-h_y(p_0)\bigr\|
\over\|y-y_0\|}=0.
\neweq\eqfinale
$$
Now, \equno\ yields
$$
\eqalign{{\bigl\|\bigl(\id-d(f_{y_0})_{p_0}\bigr)\cdot(p_y-p_0)-h_y(p_0)\bigr\|
\over\|y-y_0\|}&={\bigl\|\bigl(\id-A(y)\bigr)\cdot(p_y-p_0)-h_y(p_0)+
\bigl(A(y)-d(f_{y_0})_{p_0}\bigr)\cdot(p_y-p_0)\bigr\|
\over\|y-y_0\|}\cr
&\le{\|h_y(p_y)-h_y(p_0)\|\over\|y-y_0\|}+\bigl\|A(y)-d(f_{y_0})_{p_0}\bigr\|\,
{\|p_y-p_0\|\over\|y-y_0\|}\,.
\cr}
\neweq\eqdue
$$
Since $h_y(z)$ is holomorphic both in~$y$ and in~$z$, we have
$$
h_y(z)-h_{y_1}(z_1)=O\bigl(\|y-y_1\|,\|z-z_1\|\bigr);
$$
in particular,
$$
h_y(z)=h_y(z)-h_{y_0}(z)=O(\|y-y_0\|)
\neweq\eqtre
$$
uniformly on~$P_1$. So \equno\ implies that $p_y-p_0=O(\|y-y_0\|)$, and thus the second
summand in \eqdue\ tends to zero as $y\to y_0$.

Finally, if we apply Lemma~\Abate\ to $h_y$ we get a matrix~$B(y)$ and a constant $C>0$ such
that
$$
\|h_y(p_y)-h_y(p_0)\|\le \|B(y)\|\cdot\|p_y-p_0\|\le C\|h_y\|_{P_2}\|p_y-p_0\|
$$
when $y$ is close enough to~$y_0$, where $P_2\subset\subset P_1$ is a fixed polydisk centered
at~$p_0$. But then \eqtre\ yields
$$
\|h_y(p_y)-h_y(p_0)\|=O\bigl(\|y-y_0\|^2\bigr),
$$
and so \eqfinale\ is proved.\qedn

If $X$ is a bounded domain in~$\C^n$, then $\Hol_c(X,X)$ is an open subset of
$H^\infty(X)^n$, the complex Banach space of $n$-uples of bounded holomorphic functions
defined on~$X$. Therefore in this case $\Hol_c(X,X)$ has a natural complex structure, and we
obtain the following generalization of the main result in~[B]:

\newthm Corollary \finale: Let $D\subset\subset\C^n$ be a bounded domain. Then the Heins map
$\tau\colon\Hol_c(D,D)\to D$ is holomorphic.

\pf It follows from Theorem~\Bracci\ and [FV, Theorem~I\negthinspace I.3.10].\qedn

\setref{LLL}
\beginsection References

\art A1 M. Abate: Horospheres and iterates of holomorphic maps! Math. Z.! 198 1988 225-238

\book A2 M. Abate: Iteration theory of holomorphic maps on taut manifolds! Mediterranean
Press, Cosenza, 1989

\art A3 M. Abate: A characterization of hyperbolic manifolds!
Proc. Amer. Math. Soc.! 117 1993 789-793

\art B F. Bracci: Holomorphic dependence of the Heins map! Complex
Variables Theory and Appl.! 47 2002 1115-1119

\art C H. Cartan: Sur les r\'etractions d'une vari\'et\'e! C.R. Acad. Sci. Paris! 303 1986
715-716

\book FV T. Franzoni, E. Vesentini: Holomorphic maps and invariant distances! North Holland,
Amsterdam, 1980

\art H M.H. Heins: On the iteration of functions which are
analytic and single-valued in a given multiply-connected region!
Amer. J. Math.! 63 1941 461-480

\book He M. Herv\'e: Several complex variables. Local theory! Oxford University Press,
London, 1963

\art JK J.E. Joseph, M.H. Kwack: A generalization of a theorem of
Heins! Proc. Amer. Math. Soc.! 128 2000 1697-1701

\art R J.F. Ritt: On the conformal mapping of a region into a part of itself! Ann. of Math.!
22 1920 157-160

\art Ro H. Rossi: Vector fields on analytic spaces! Ann. of Math.! 78 1963 455-467

\art W R. Wavre: Sur la r\'eduction des domaines par une substitution \`a $m$ variables
complexes et l'existence d'un seul point invariant! Enseign. Math.! 25 1926 218-234

\bigbreak
\leftline{\eightrm Marco Abate}
\leftline{\eightrm Dipartimento di Matematica, Universit\`a di Pisa}
\leftline{\eightrm Largo Pontecorvo 5, 56127 Pisa, Italy}
\leftline{\eightrm E-mail: abate@dm.unipi.it}
\medskip
\leftline{\eightrm Filippo Bracci}
\leftline{\eightrm Dipartimento di Matematica, Universit\`a di Roma
``Tor Vergata''}
\leftline{\eightrm Via della Ricerca Scientifica, 00133 Roma, Italy}
\leftline{\eightrm E-mail: fbracci@mat.uniroma2.it}
\bye